\documentclass[10pt]{article}
\usepackage{amssymb, amsmath, url, graphicx, setspace, geometry}
\usepackage{calrsfs, xcolor}
\usepackage{wasysym}
\usepackage{listings}

\def\3{\subset }
\def\4{\subseteq }
\def\<{\left<}
\def\>{\right>}

\def\bit{\begin{itemize}}
\def\eit{\end{itemize}}
\def\3{\subset }
\def\4{\subseteq }

\def\0{\leqno}

\def\barr{\begin{array}}
\def\earr{\end{array}}

\def\Z{{\rlap{$\kern2pt{\rm Z}$}{\rm Z}\,}}
\def\bld#1#2{{\buildrel{#1}\over{#2}}}
\def\st#1#2{{\mathrel{\mathop{#2}\limits_{#1}}{}\!}}
\def\stb#1#2#3{{\st{{#1}}{\bld{{#2}}{#3}}{}\!}}
\def\xmare#1#2{\stb{#1}{#2}{\mbox{\Large$\times$}}}

\title{\bf On the number of conjugacy classes of subgroups of a finite group}
\author{Mihai-Silviu Lazorec and Marius T\u arn\u auceanu}
\date{August 6, 2026}

\begin{document}

\maketitle

\begin{abstract}
Let $k'(G)$ and $L(G)$ be the number of conjugacy classes of subgroups and the subgroup lattice of a finite group $G$, respectively. Our objective is to study some aspects related to  the ratios $d'(G)=\frac{k'(G)}{|L(G)|}$ and $d^*(G)=\min\{ d'(S) \mid S\text{ is a section of }G\}$ which measure how close $G$ is to being a Dedekind group.  We prove that the set containing the values $d'(G)$, as $G$ ranges over the class of nilpotent groups, is dense in $[0, 1]$. A nilpotency criterion is  obtained by proving that if $d^*(G)>\frac{2}{3}$, then $G$ is nilpotent and   information on its structure is given. We also show that if $d^*(G)>\frac{4}{5}$, then $G$ is an Iwasawa group. Finally, we deduce a result which ensures that a $p$-group of order $p^n$ ($n\geq 3$) is a Dedekind group. This last result can be extended to the class of nilpotent groups and it also highlights the second maximum value of $d^*$ on the class of $p$-groups of order $p^n$.   
\end{abstract}



\section{Introduction}

All groups considered in this paper are finite and $n\geq 1$ denotes an integer. Let $G$ be a group and $p$ be a prime number. We denote by $k(G), k'(G), \nu(G), N(G)$ and $L(G)$ the number of conjugacy classes, the number of conjugacy classes of subgroups, the number of conjugacy classes of non-normal subgroups, the normal subgroup lattice and the subgroup lattice of $G$, respectively. The cyclic group of order $n$ is denoted by $C_n$, the dihedral group $D_{2n}$ ($n\geq 3$) has the following structure 
$$D_{2n}=\langle x, y \mid x^n=y^2=1, yx=x^{n-1}y\rangle,$$
$M_{p^n}$ ($n\geq 4$ if $p=2$; $n\geq 3$ if $p$ is odd) is a modular $p$-group (i.e. a $p$-group whose subgroup lattice is modular) given by
$$M_{p^n}=\langle x, y \mid x^{p^{n-1}}=y^p=1,  yx=x^{p^{n-2}+1}y\rangle,$$
while $He_p$ is the Heisenberg group modulo $p$ ($p$ is odd) and its representation is
$$He_p=\langle x, y, z \ | \ x^p=y^p=z^p=1, [x, z]=[y, z]=1, [x, y]=z\rangle.$$
We also recall that an Iwasawa group is a modular nilpotent group.

A way to measure how close $G$ is to being abelian is given by its commutativity degree $$d(G)=\frac{|\{ (x, y)\in G\times G \mid xy=yx\}|}{|G|^2}.$$
Theorem IV of \cite{6} states that the numerator is equal to $k(G)|G|$, so the commutativity degree of $G$ can be also expressed as
$$d(G)=\frac{k(G)}{|G|}.$$
In this paper, we study another ratio that is defined similarly, the main difference being that we switch from working with the elements of $G$ to operating with its subgroups. Hence,  we are going to investigate the ratio
$$d'(G)=\frac{k'(G)}{|L(G)|}.$$
This quantity measures how close is $G$ from being a Dedekind group (i.e. all its subgroups are normal). It is clear that all abelian groups are Dedekind, while the structure of the non-abelian Dedekind groups, which are  called Hamiltonian groups, is also known (see statement 5.3.7 of \cite{12}). 

Some properties of $d'$ are straightforward. For any group $G$, we have $0<d'(G)\leq 1$, and the equality $d'(G)=1$ holds if and only if $G$ is a Dedekind group. If $G_1$ and $G_2$ are two groups such that $G_1\cong G_2$, then $d'(G_1)=d'(G_2)$. The converse does not hold even for non-Dedekind groups. For instance, if $G_1\cong C_3\rtimes C_8$ (SmallGroup(24,1)) and $G_2\cong C_3\times D_8$ (SmallGroup(24,10)), we have $d'(G_1)=d'(G_2)=\frac{4}{5}$ and $G_1\not\cong G_2$. Finally, if $I$ is a finite set with $|I|\geq 2$ and $(G_i)_{i\in I}$ is a family of groups of coprime orders, then $$d'\bigg(\xmare{i\in I}{} G_i\bigg)=\prod_{i\in I}d'(G_i),$$ so $d'$ is  multiplicative.

Our main results are listed below and their proofs are included in the following two sections. In Section 2, we explicitly compute $d'(G)$ with $G$ belonging to specific classes of groups. Given an integer $a\geq 1$, we show that all ratios of the form $\frac{a}{a+1}$ are attained by $d'$. The main result of this section is the following one:\\

\textbf{Theorem A.} \textit{Let $\mathcal{N}$ be the class of nilpotent groups. Then the set $$\{ d'(G) \mid G\in\mathcal{N}\}$$ is dense in $[0, 1]$.}\\

In Section 3, we justify that there is no constant $c\in (0, 1)$ such that if $d'(G)>c$, then $G$ is a nilpotent/Iwasawa/Dedekind group. However, a nilpotency criterion and a criterion for $G$ to be an Iwasawa group can be obtained by replacing $d'(G)$ with $d^*(G)$, where 
$$d^*(G)=\min\{ d'(S) \mid S\text{ is a section of }G\}.$$
More exactly, we justify that the following results hold:\\

\textbf{Theorem B.} \textit{Let $G$ be a group. If $d^*(G)>\frac{2}{3}=d^*(S_3)=d'(S_3)$, then $G$ is a nilpotent group. Moreover, if $2 \nmid |G|$, then $G$ is the direct product of some modular $p$-groups.}\\

\textbf{Theorem C.} \textit{Let $G$ be a group. If $d^*(G)>\frac{4}{5}=d^*(D_8)=d'(D_8)$, then $G$ is an Iwasawa group.}\\

By using the classification of minimal non-Dedekind $p$-groups, we also determine a criterion which guarantees that $p$-groups of fixed orders are Dedekind groups. This result can be extended to the class of nilpotent groups and it is stated below:\\

\textbf{Theorem D.} \textit{Let $n\geq 3$ be an integer and $G$ be a $p$-group of order $p^n$. 
\begin{itemize}
\item[i)] Assume that $p=2$ and $n=3$. If $d^*(G)>\frac{4}{5}=d^*(D_8)=d'(D_8)$, then $G$ is a Dedekind group;
\item[ii)] Assume that $n\geq 4$ if $p=2$, and $n\geq 3$ if $p$ is odd. If $d^*(G)>\frac{(n-2)(p+1)+4}{(n-1)(p+1)+2}=d^*(M_{p^n})=d'(M_{p^n})$, then $G$ is a Dedekind group.
\end{itemize}}

The paper concludes with a list of open problems.
  
\section{On the values of \textit{d'}}

As a preliminary result we recall the main theorem of \cite{3}. It outlines the  classification of groups having one conjugacy class of non-normal subgroups. These groups are close to being Dedekind groups and this aspect will be helpful in obtaining some of our results.\\

\textbf{Lemma 2.1.} \textit{Let $G$ be a group. Then $\nu(G)=1$ if and only if $G$ is isomorphic to one of the following groups:
\begin{itemize}
\item[i)] $M_{p^n}$;
\item[ii)] a non-abelian split extension $G_{p, q, n}=N\rtimes P$ where $N\cong C_p$, $P\cong C_{q^{n-1}}, [N, \Phi(P)]=1$, $p, q$ are primes such that $q\mid p-1$ and $n\geq 2$.
\end{itemize}}

In what follows we obtain explicit formulas for computing $d'(G)$ when $G$ is a group with $\nu(G)=1$, a dihedral 2-group or a Heisenberg group modulo $p$.\\

\textbf{Proposition 2.2.} \textit{The following results hold:
\begin{itemize}
\item[i)] $d'(M_{p^n})=\frac{(n-2)(p+1)+4}{(n-1)(p+1)+2};$
\item[ii)] $d'(G_{p,q,n})=\frac{2n}{2n+p-1};$
\item[iii)] $d'(D_{2^n})=\frac{3n-1}{2^n+n-1};$
\item[iv)] $d'(He_p)=\frac{2p+5}{p^2+2p+4}.$
\end{itemize}
Moreover, the sequences $(d'(M_{p^n}))_{n\geq 3}$ ($p$ is odd), $(d'(M_{2^n}))_{n\geq 4}$, $(d'(G_{p, q, n}))_{n\geq 2}$ are strictly increasing, while the sequences $(d'(D_{2^n}))_{n\geq 3}$, $(d'(He_p))_{p\geq 3}$ are strictly decreasing.}\\

\textbf{Proof.} Note that $k'(G)=|N(G)|+\nu(G)$ for any group $G$. According to Lemma 3.3 of \cite{15}, we have
$$|L(M_{p^n})|=(n-1)(p+1)+2 \text{ and } |L(D_{2^n})|=2^n+n-1.$$

\textit{i)} As stated in Lemma 3.5 of \cite{15}, we have $$|N(M_{p^n})|=(n-2)(p+1)+3.$$ Since $\nu(M_{p^n})=1$, the conclusion follows.

\textit{ii)} The subgroup $P\cong C_{q^{n-1}}$ is the representative of the unique conjugacy class of non-normal subgroups of $G_{p, q, n}$. This class contains $p$ conjugates. Besides $N$ and its improper subgroups, $G_{p, q, n}$ has one normal subgroup of order $pq^i$ and one normal subgroup of order $q^i$ for each $i\in \{ 1, 2, \ldots, n-2\}$. Therefore, 
$$k'(G_{p, q, n})=2n \text{ and } |L(G_{p, q, n})|=2n+p-1,$$
which leads to the desired result.  

\textit{iii)} Theorem 3.3 of \cite{5} outlines the conjugacy classes of subgroups of any dihedral group. According to this result, $D_{2^n}$ has 3 conjugacy classes of subgroups of order $2^i$ for each $i\in\{ 1,2,\ldots, n-1\}$. Hence, $$k'(D_{2^n})=3n-1,$$
and the conclusion follows.

\textit{iv)} The proper non-trivial normal subgroups of $He_p$ are $Z(He_p)\cong C_p$ and its $p+1$ maximal subgroups isomorphic to $C_p^2$. There are also $p+1$ non-trivial conjugacy classes of subgroups, each being formed of $p$ conjugates isomorphic to $C_p$. We conclude that
$$k'(He_p)=2p+5 \text{ and } |L(He_p)|=p^2+2p+4,$$ 
so we easily obtain the value of $d'(He_p)$.

Regarding the second part of this proposition, we only justify the monotonicity of the first sequence. Let $p$ be a fixed odd prime and  $f:[3, \infty)\longrightarrow \mathbb{R}$ be a function given by $$f(x)=\frac{(x-2)(p+1)+4}{(x-1)(p+1)+2}, \ \forall \ x\in [3, \infty).$$
We have
$$f'(x)=\frac{p^2-1}{[(x-1)(p+1)+2]^2}>0, \ \forall \ x\in (3, \infty),$$
so $f$ is strictly increasing on $[3, \infty)$. Hence, the sequence  $(f(n))_{n\geq 3}=(d'(M_{p^n}))_{n\geq 3}$ is strictly increasing. 
\hfill\rule{1,5mm}{1,5mm}\\

The following result is a consequence of our  previous proposition. It shows that $d'(G)$ can be made arbitrarily large or small as $|G|$ increases.\\

\textbf{Corollary 2.3.} \textit{The following results hold:
$$\displaystyle\lim_{n\to\infty}d'(M_{p^n})=\displaystyle\lim_{n\to\infty}d'(G_{p, q, n})=1 \text{ and } \displaystyle\lim_{n\to\infty}d'(D_{2^n})=\displaystyle\lim_{p\to\infty}d'(He_p)=0.$$}

It is clear that the values of $d'$ belong to $S=[0, 1]\cap \mathbb{Q}$. We proceed our study by showing that specific values from $S$ are attained by $d'$. Then, in the proof of Theorem A, we justify that each value of $S$ is an adherent point of the set
$\{ d'(G) \mid G\in\mathcal{N}\},$
where $\mathcal{N}$ is the class of nilpotent groups. The classification given by Lemma 2.1 and the numerical results outlined in Proposition 2.2 play a significant role for both purposes.\\

\textbf{Corollary 2.4.} \textit{Let $a\geq 1$ be an integer. Then there is a group $G$ such that $d'(G)=\frac{a}{a+1}.$}\\

\textbf{Proof.} If $a=1$, we take $G=G_{5, 2, 2}\cong D_{10}$ and we have $d'(G)=\frac{1}{2}$. Suppose that $a\geq 2$. Then, we can choose $G=G_{3, 2, a}$ and we get $d'(G)=\frac{a}{a+1}.$
\hfill\rule{1,5mm}{1,5mm}\\ 

\textbf{Proof of Theorem A.} Let $x\in S$. We aim to show that there is a sequence of groups $(G_n)_{n\geq n_0}\subseteq\mathcal{N}$ such that $$\displaystyle\lim_{n\to\infty}d'(G_n)=x.$$  

If $x=0$ or $x=1$, we are done by Corollary 2.3. Suppose that $x=\frac{a}{b}$ where $1\leq a<b$. Also, let $(p_k)_{k\geq 1}$ be the strictly increasing sequence of odd prime numbers. We construct $b-a$ strictly increasing and disjoint subsequences $(p_{1_n})_{n\geq 1}, (p_{2_n})_{n\geq 1}, \ldots, (p_{{(b-a)}_n})_{n\geq 1}$ of $(p_k)_{k\geq 1}$. Also, for each $i\in\{ 1, 2, \ldots, b-a \}$, consider the sequence of groups $(G_{i_n})_{n\geq 1}$ where $$G_{i_n}=M_{{p_{i_n}}^{a+i+1}}.$$ 
Note that each of these $b-a$ sequences of groups is well-defined since we work only with odd prime numbers and $a+i+1\geq 3$. Finally, we define the sequence $(G_n)_{n\geq 1}\subseteq\mathcal{N}$ whose general term is 
$$G_n=\xmare{i=1}{b-a}G_{i_n}.$$
The factors of the above direct product are of coprime orders since the sequences $(p_{i_n})_{n\geq 1}$ are disjoint. Hence, we are able to apply the multiplicativity property of $d'$ as follows
$$d'(G_n)=d'\bigg(\xmare{i=1}{b-a}G_{i_n}\bigg )=\prod\limits_{i=1}^{b-a}d'(G_{i_n})=\prod\limits_{i=1}^{b-a}d'(M_{{p_{i_n}}^{a+i+1}})=\prod\limits_{i=1}^{b-a}\frac{(a+i-1)(p_{i_n}+1)+4}{(a+i)(p_{i_n}+1)+2}.$$
As $n$ tends to infinity, we obtain
$$\displaystyle\lim_{n\to\infty}d'(G_n)=\prod\limits_{i=1}^{b-a}\displaystyle\lim_{n\to\infty}\frac{(a+i-1)(p_{i_n}+1)+4}{(a+i)(p_{i_n}+1)+2}=\prod\limits_{i=1}^{b-a}\frac{a+i-1}{a+i}=\frac{a}{b}=x.$$
Therefore, we conclude that each value of $S$ is an adherent point of the set $\{ d'(G) \ | \ G\in\mathcal{N}\}$. This implies that
$$S\subseteq \overline{\{ d'(G) \mid G\in\mathcal{N}\}},$$
so $$[0, 1]=\overline{S}\subseteq \overline{\{ d'(G) \mid  G\in\mathcal{N}\}}.$$
Since the reverse inclusion is obvious, the proof is complete.
\hfill\rule{1,5mm}{1,5mm}\\ 

Theorem A can be easily  extended to any class of groups containing $\mathcal{N}$ as follows.\\

\textbf{Corollary 2.5.} \textit{Let $\mathcal{C}$ be a class of groups such that $\mathcal{N}\subseteq \mathcal{C}$. Then $\{ d'(G) \mid G\in\mathcal{C}\}$ is dense in $[0, 1].$} 

\section{Criteria ensuring that a finite group is nilpotent / Iwasawa / Dedekind}

Let $G$ be a group. Given a group theoretic property $P$, $G$ is said to be a minimal non-$P$ group if $G$ does not satisfy $P$ and all its proper subgroups do. The following lemma recalls the structure and some properties of the minimal non-nilpotent groups. These are also known as Schmidt groups. The listed results are outlined in \cite{13} and \cite{11}.\\

\textbf{Lemma 3.1.} \textit{Let $G$ be a  Schmidt group. Then the following statements hold:
\begin{itemize}
\item[i)] $G=P\rtimes Q$ where $P$ is a normal Sylow $p$-subgroup, $Q$ is a cyclic Sylow $q$-subgroup of $G$ and $p, q$ are distinct primes;
\item[ii)] $Z(G)=\Phi(G)=\Phi(P)\times\Phi(Q)$;
\item[iii)] $\frac{P}{\Phi(P)}\cong C_p^r$, where $r$ is the multiplicative order of $p$ modulo $q$;
\item[iv)] If $N$ is a proper normal subgroup of $G$, then $N$ does not contain $Q$ and either $P\subseteq N$ or $N\subseteq Z(G)$. 
\end{itemize}}

Let $c\in (0, 1)$ and $G$ be a group. A condition such as $d'(G)>c$ cannot guarantee that $G$ is a nilpotent/Dedekind group. Indeed, according to Corollary 2.3 for a sufficiently large value of $n$, say $n_c$, we would get $d'(M_{p^{n_c}})>c$ and $d'(G_{p, q, n_c})>c$, but $M_{p^{n_c}}$ is a non-Dedekind group, while $G_{p, q, n_c}$ is non-nilpotent.   Theorem B highlights a nilpotency criterion that is obtained by replacing the condition $d'(G)>c$ with $d^*(G)>c$, where
$$d^*(G)=\min\{ d'(S) \mid S\text{ is a section of } G\}.$$
A useful property of $d^*$ is that it behaves like a decreasing function defined on $L(G)$, i.e. if $H, K\in L(G)$ and $H\subseteq K$, then $d^*(H)\geq d^*(K).$ Also, $d^*$ is multiplicative. Indeed, if $I$ is a finite set with $|I|\geq 2$ and $(G_i)_{i\in I}$ is a family of groups of coprime orders, then all sections of $\xmare{i\in I}{}G_i$ are of the form $\xmare{i\in I}{}S_i$, where $S_i$ is a section of $G_i$ for any $i\in I$. Then
$$d^*\bigg(\xmare{i\in I}{} G_i\bigg)=\prod_{i\in I}d^*(G_i).$$

Before proving Theorem B, we outline and justify some modularity criteria for $p$-groups.\\

\textbf{Proposition 3.2.} \textit{Let $G$ be a $p$-group.
\begin{itemize}
\item[i)] If $d^*(G)>\frac{4}{5}=d^*(D_8)=d'(D_8)$, then $G$ is a modular $p$-group;
\item[ii)] If $p$ is odd and $d^*(G)>\frac{11}{19}=d^*(He_3)=d'(He_3)$, then $G$ is a modular $p$-group.
\end{itemize}}

\textbf{Proof.} \textit{i)} Suppose that $G$ is a non-modular $p$-group. According to Lemma 2.3.3 of \cite{14}, $G$ has a section $S$ isomorphic to:
\begin{itemize}
\item[--] $D_8$, if $p=2$;
\item[--] $He_p$, if $p$ is odd.
\end{itemize}
We also know that $d^*(G)\leq d'(S)$. We use the results stated in Proposition 2.2. If $p=2$, we get $$d^*(G)\leq \frac{4}{5}=d'(D_8),$$
a contradiction. If $p$ is odd, then
$$d^*(G)\leq d'(He_p)\leq \frac{11}{19}=d'(He_3)<\frac{4}{5},$$
a contradiction. 

The above reasoning also justifies statement  \textit{ii)}, so the proof is complete.    
\hfill\rule{1,5mm}{1,5mm}\\ 

We proceed by proving Theorem B.\\

\textbf{Proof of Theorem B.} Suppose that $G$ is a counterexample of minimal order. Then $G$ is non-nilpotent. Let $H\in L(G)$ such that $H\neq G$. Then
$$d^*(H)\geq d^*(G)>\frac{2}{3},$$
so $H$ is nilpotent. Consequently, $G$ is a Schmidt group. 

We consider the section $S=\frac{G}{Z(G)}$ of $G$. By using the information introduced in Lemma 3.1, we deduce that $|S|=p^rq$ and, by the correspondence theorem, $S$ has a normal subgroup $$H=\frac{P\times \Phi(Q)}{Z(G)}\cong C_p^r.$$ On the other hand, a subgroup $K\cong C_q$ of $S$ cannot be normal. Indeed, if $K\in N(S)$, then $G$ would have a proper normal subgroup $N$ with $|N|=q\cdot |Z(G)|=|\Phi(P)|\cdot |Q|$. By Lemma 3.1 \textit{iv)}, we would also have that $|P|$ divides $|N|$ and this leads to a contradiction. Hence, we conclude that
$$S=H\rtimes K\cong C_p^r\rtimes C_q$$
and $K$ acts faithfully on $H$. 
Moreover, all proper subgroups of $S$ are nilpotent, so $S$ is also a Schmidt group. By Lemma 3.1 \textit{ii)}, we have $Z(S)=\{ 1\}$.  

Note that $S$ cannot have a subgroup $L_i$ of order $p^iq$ with $1\leq i<r$. If such a subgroup would exist, then it would be nilpotent, so we can assume that there is a subgroup $H_i\cong C_p^i$ of $H$ such that $L_i=H_i\times K\cong C_p^i \times C_q$.
Then $K$ centralizes $H_i$. Since $H$ also centralizes $H_i$, it follows that $S\subseteq C_S(H_i)$. Then $H_i\subseteq Z(S)=\{1\}$, a contradiction.

We denote $|L(C_p^r)|$ by $a_{p, r}$. Taking into account that $S$ has $p^r$ conjugate subgroups isomorphic to $C_q$ and there are no    subgroups of order $p^iq$ with $1\leq i<r$, we obtain
$$|L(S)|=a_{p, r}+p^r+1.$$  
Let $i\in \{1, 2, \ldots, r-1\}$ and let $H_i$ be a subgroup of $S$ with $|H_i|=p^i$. By our previous discussion on the non-existence of subgroups of order $p^iq$, we deduce that $H_i$ is not normal in $S$. Also, $H$ is an abelian maximal subgroup of $S$ and $H_i\subset H$, so the size of the conjugacy class of $H_i$ is
$$[S:N_S(H_i)]=[S:H]=q.$$
It is known that the number of subgroups of $S$ that are isomorphic to $H_i$ is given by the Gaussian coefficient $$\big[\substack{r \\ i}\big]_p=\frac{(p^r-1)(p^{r-1}-1)\cdots (p^{r-i+1}-1)}{(p^i-1)(p^{i-1}-1)\cdots (p-1)}$$ (see, for example, Exercise 1.74 of \cite{1}). Hence, these subgroups are partitioned into $\frac{1}{q}\big[\substack{r \\ i}\big]_p$ conjugacy classes. Besides them, $S$ has 3 trivial conjugacy classes of subgroups whose representatives are $1$, $H$ and $S$, respectively, and the conjugacy class of $K$. Therefore,
$$k'(S)=\frac{1}{q}\sum\limits_{i=1}^{r-1}\big[\substack{r \\ i}\big]_p+4=\frac{a_{p, r}+4q-2}{q}.$$ 
Consequently, we have
$$d'(S)=\frac{a_{p, r}+4q-2}{q(a_{p, r}+p^r+1)}.$$ 

The following step is to show that $d'(S)\leq\frac{2}{3}$. This would contradict that $d^*(G)>\frac{2}{3}$ and the first part of our proof would be complete. Note that
$$
d'(S)\leq\frac{2}{3}\Longleftrightarrow (2q-3)a_{p,r}+2q(p^r-5)+6\geq 0
$$
Since $q\geq 2$, in most cases it suffices to justify that
\begin{equation}\label{r1}
p^r-5\geq 0.
\end{equation}
This is true for $p\geq 5$, so there are two cases left to discuss: $p=2$ and $p=3$. We recall that $q\mid p^r-1$. If we assume that $p=2$, then $r\geq 2$. If $r\geq 3$, then (\ref{r1}) holds, while if $r=2$, then $q=3$ and $S\cong C_2^2\rtimes C_3\cong A_4$, so $d'(S)=\frac{1}{2}<\frac{2}{3}$. Finally, suppose that $p=3$. If $r\geq 2$, then (\ref{r1}) holds, while if $r=1$, then $q=2$ and $S\cong C_3\rtimes C_2\cong S_3$, so $d'(S)=\frac{2}{3}$. Therefore, our initial assumption is false, so  $G$ is a nilpotent group.

We now prove the second part of the statement of Theorem B. Since $G$ is nilpotent and $2\nmid |G|$, it follows that $G\cong\xmare{i=1}{k}G_i$, where $G_1, G_2, \ldots, G_k$ are the Sylow subgroups of $G$, each of them being of odd order. We have
$$d^*(G_i)\geq d^*(G)>\frac{2}{3}>\frac{11}{19}, \ \forall \ i\in\{ 1, 2, \ldots, k\}.$$
Then each $G_i$ is a modular group due to Proposition 3.2 \textit{ii)}. Hence, $G$ is a direct product of modular $p$-groups.
\hfill\rule{1,5mm}{1,5mm}\\ 

Regarding the second part of the statement of Theorem B, we mention that the structure of modular $p$-groups is determined by Iwasawa in \cite{9} (see also Theorem 2.3.1 of \cite{14}). More exactly, $P$ is a modular $p$-group if and only if one of the following holds: 
\begin{itemize}
\item[--] $P\cong Q_8\times C_2^n$ ($n\geq 0$);
\item[--] $P$ contains an abelian normal subgroup $A$ such that $\frac{P}{A}$ is cyclic; in addition, there is $g\in P$ and a positive integer $n$ such that $P=A\langle g\rangle$ and $ga=a^{p^n+1}g$ for all $a\in A$, with $n\geq 2$ if $p=2$.
\end{itemize}
Also, according to Exercise 3 from Section 2.4 of \cite{14}, the fact that $G$ is a direct product of modular $p$-groups is equivalent to  any of the following:
\begin{itemize}
\item[--] $G$ is an Iwasawa group;
\item[--] all subgroups of $G$ are permutable.
\end{itemize} 
These additional details are useful for proving Theorem C.\\

\textbf{Proof of Theorem C.} Since $$d^*(G)>\frac{4}{5}>\frac{2}{3},$$ $G$ is a nilpotent group according to Theorem B. Then $G$ is the direct product of its Sylow subgroups. Let $P$ be one of them. We have $$d^*(P)\geq d^*(G)>\frac{4}{5},$$ so $P$ is a modular $p$-group by Proposition 3.2 \textit{i)}. Consequently, $G$ is an Iwasawa group, as desired.   
\hfill\rule{1,5mm}{1,5mm}\\ 

We mention that the lower bounds that appear in the statements of Proposition 3.2 and Theorems B and C are the best possible ones since $D_8$ and $He_3$ are non-modular $p$-groups, while $S_3$ is a non-nilpotent group.

A classification of minimal non-Dedekind groups is obtained in \cite{4} (see also Lemma 3.1 of \cite{8}). The following result  illustrates the structure of minimal non-Dedekind $p$-groups only. The first two types of groups are also minimal non-abelian $p$-groups by Exercise 8a from Section 1 of \cite{2}. The structure of their centers  is described in Lemma 2.2 of \cite{7} and is also recalled below.\\

\textbf{Lemma 3.3.} \textit{Let $G$ be a minimal non-Dedekind $p$-group. Then $G$ is isomorphic to one of the following groups:
\begin{itemize}
\item[i)] $H_{p, s, t}=\langle x, y, z \mid x^{p^s}=y^{p^t}=z^p=1, [x, z]=[y, z]=1, [x, y]=z\rangle,$ where $s\geq t\geq 1$ and, if $p=2$, then $n=s+t\geq 3$; moreover, $Z(H_{p,s,t})=\langle x^p \rangle\times\langle y^p\rangle\times\langle z\rangle$;
\item[ii)] $K_{p, s, t}=\langle x, y\mid x^{p^s}=y^{p^t}=1, yx=x^{p^{s-1}+1}y\rangle$, where $s\geq 2, t\geq 1$ and, if $p=2$, then $n=s+t\geq 4$; moreover $Z(K_{p, s, t})=\langle x^p \rangle\times\langle y^p\rangle;$
\item[iii)] $Q_{16}$.
\end{itemize}
}

A group $G$ is a q-self dual group if every quotient of $G$ is isomorphic to a subgroup of $G$. Theorems 3 and 4 of \cite{10} list the $p$-groups $G$ all of whose subgroups are $q$-self-dual. For odd primes, the obtained classification is valid under the assumption that $\Omega_1(G)$ is abelian. The following preliminary result is a consequence of these theorems.\\

\textbf{Lemma 3.4.} \textit{$M_{p^n}$ is a q-self dual group.}\\

The following result shows that the values of $d'$ and $d^*$ are equal for some $p$-groups.\\

\textbf{Proposition 3.5.} \textit{The following equalities hold:
\begin{itemize}
\item[i)] $d'(M_{p^n})=d^*(M_{p^n})$;
\item[ii)] $d'(He_p)=d^*(He_p)$.
\end{itemize}}

\textbf{Proof.} \textit{i)} Let $G=M_{p^n}\cong K_{p, n-1, 1}$. Then $G$ is a minimal non-Dedekind group by Lemma 3.3. Let $H, K\in L(G)$ such that $K\in N(H)$ and $H\neq G$. Then $H$ is a Dedekind group and this property is inherited by the section $\frac{H}{K}$. Hence, we deduce that
$$d^*(G)=\min\bigg\{ d'\bigg(\frac{G}{N}\bigg)\mid N\in N(G)\bigg\}.$$
Since $G$ is a q-self dual group and all its proper subgroups are abelian, $\frac{G}{1}\cong G$ is its only non-Dedekind section. Therefore, $d^*(G)=d'(G).$

\textit{ii)} The only non-Dedekind section of $G\cong He_p$ is isomorphic to $G$, so the conclusion follows.  
\hfill\rule{1,5mm}{1,5mm}\\ 

Let $c\in(0, 1)$. The idea of replacing the hypothesis $d'(G)>c$ with $d^*(G)>c$ in order to obtain a criterion for $G$ to be a Dedekind group does not work. Indeed, due to Proposition 3.5, the two conditions are equivalent for $G\cong M_{p^n}$ and we explained that there is a non-Dedekind group $M_{p^{n_c}}$ such that $d'(M_{p^{n_c}})=d^*(M_{p^{n_c}})>c$. However, we are able to show that such a criterion can be obtained if we work with $p$-groups of fixed orders. Theorem D illustrates this result and our next objective is to prove it. We justify a preliminary result first.\\

\textbf{Lemma 3.6.} \textit{Let $p$ be a prime number. The following results hold:
\begin{itemize}
\item[i)] If $p=2$ and $n=s+t\geq 3$, then $H_{2, s, t}$ has a section isomorphic to $D_8$;
\item[ii)] If $p$ is odd and $n=s+t\geq 2$, then $H_{p, s,t}$ has a section isomorphic to $He_p$;
\item[iii)] If $p=2, n=s+t\geq 5$ and $K_{2, s, t}\not\cong M_{2^n}$, then $K_{2, s, t}$ has a section isomorphic to $M_{2^k}$, where $k<n$;
\item[iv)] If $p$ is odd, $n=s+t\geq 4$ and $K_{p,s,t}\not\cong M_{p^n}$, then $K_{p, s, t}$ has a section isomorphic to $M_{p^k}$, where $k<n$. 
\end{itemize}}

\textbf{Proof.} By Lemma 3.3, we know that 
$$N=\langle x^p\rangle\times \langle y^p\rangle\cong C_{p^{s-1}}\times C_{p^{t-1}}  \text{ and } M=\langle y^{p^{t-1}}\rangle\cong C_p$$ 
are central subgroups of $H_{p, s, t}$ and $K_{p, s, t}$, respectively. Hence, they are  also normal subgroups.  

\textit{i)} \& \textit{ii)} Let $G=H_{p, s, t}$, $\bar{x}=xN, \bar{y}=yN$ and $\bar{z}=zN$. Then 
$$S=\frac{G}{N}=\langle \bar{x}, \bar{y}, \bar{z} \mid \bar{x}^p=\bar{y}^p=\bar{z}^p=N, [\bar{x}, \bar{z}]=[\bar{y}, \bar{z}]=N, [\bar{x}, \bar{y}]=\bar{z} \rangle$$
is a non-abelian section of order $p^3$ of $G$. If $p=2$, it is easy to check that $S$ has 5 elements of order 2, so $S\cong D_8$. Otherwise, $p$ is odd and $S\cong He_p$. 

\textit{iii)} We proceed by induction on $n$. For the initial step, if $n=5$, then one can use GAP \cite{16}\footnote{See Appendix \ref{gap} for the GAP code used to verify these results.} to check that $K_{2, 3, 2}$ (SmallGroup(32, 4)) and $K_{2, 2, 3}$ (SmallGroup(32, 12)) have a section isomorphic to $M_{16}$.

Let $n\geq 6$ and suppose that the statement holds for any $m$ such that $5\leq m<n$. Let $G=K_{2, s, t}\not \cong M_{2^n}$. Then $t\geq 2$. We have
$$\frac{G}{M}\cong K_{2, s, t-1}.$$
If $K_{2, s, t-1}\cong M_{2^{n-1}}$, then we are done. Otherwise, by the inductive hypothesis, it follows that $\frac{G}{M}$ has a section $$S=\frac{\frac{H}{M}}{\frac{K}{M}}\cong \frac{H}{K}\cong M_{2^k},$$ where $H, K\in L(G)$ and $K\in N(H)$. Then $G$ also has a section isomorphic to $M_{2^k}$ and the proof is complete. 

\textit{iv)} This is also done by induction on $n$. If $n=4$, then $G=K_{p, 2, 2}$ has a section
$$\frac{G}{M}\cong K_{p, 2, 1}\cong M_{p^3}.$$
Thus, the initial step is complete. The inductive step follows a similar argument as the one done in the proof of item \textit{iii)}, so we omit it.  
\hfill\rule{1,5mm}{1,5mm}\\ 

We continue with proving Theorem D.\\

\textbf{Proof of Theorem D.} The result outlined in item \textit{i)} is clear since $D_8$ is the only non-Dedekind group of order 8.

Regarding item \textit{ii)}, let $G$ be a $p$-group such that 
\begin{equation}\label{r2}
d^*(G)>d^*(M_{p^n}).
\end{equation}   
We proceed by induction on $n$. For the initial step, we work with groups of order 16 or $p^3$, where $p$ is odd. By checking the list of groups of order 16, we deduce that (\ref{r2}) holds if and only if $G$ is a Dedekind group. In addition, the only non-Dedekind groups of order $p^3$ are $M_{p^3}$ and  $He_p$. By Propositions 2.2 and 3.5, we have
$$d^*(He_p)=d'(He_p)=\frac{2p+5}{p^2+2p+4}<\frac{p+5}{2p+4}=d'(M_{p^3})=d^*(M_{p^3}).$$
Therefore, (\ref{r2}) holds again if and only if $G$ is a Dedekind group. 

Let $n\geq 5$ when $p=2$, and $n\geq 4$ when $p$ is odd. Assume that the statement holds for any $m$ such that $m<n$. If $p=2$, we may assume that $m\geq 4$, whereas if $p$ is odd, then $m\geq 3$. Let $G$ be a $p$-group of order $p^n$ such that (\ref{r2}) holds. Based on the parity of $p$, by using Propositions 2.2 and 3.5 again, we have
\begin{equation*}
d^*(G)>d^*(M_{2^n})\geq d^*(M_{32})=\frac{13}{14} \text{ or } d^*(G)>d^*(M_{p^n})\geq d^*(M_{p^4})=\frac{2p+6}{3p+5}.
\end{equation*}

For the sake of contradiction, suppose that $G$ is a non-Dedekind group. Let $H$ be a maximal subgroup of $G$. By using Propositions  2.2 and 3.5, we have
$$d^*(H)\geq d^*(G)>d^*(M_{p^n})>d^*(M_{p^{n-1}}).$$
Therefore, $H$ is a Dedekind group by the inductive hypothesis. Consequently, $G$ is a minimal non-Dedekind group. By Lemma 3.3, it follows that $G\cong H_{p, s, t}$ or $G\cong K_{p, s, t}$.   

If $G\cong H_{p, s, t}$, then, according to Lemma 3.6, $G$ has a section isomorphic to $D_8$ or $He_p$. In the first case, we get $d^*(G)\leq d'(D_8)=\frac{4}{5}$, contradicting $d^*(G)>\frac{13}{14}$. In the second case, we obtain $d^*(G)\leq d'(He_p)=\frac{2p+5}{p^2+2p+4}$, contradicting $d^*(G)>\frac{2p+6}{3p+5}.$

Suppose that $G\cong K_{p, s, t}$. By (\ref{r2}), we know that $G\not\cong M_{p^n}$. Hence, we can use Lemma 3.6 to conclude that $G$ has a section isomorphic to $M_{p^k}$, where $k<n$. It follows that
$$d^*(G)\leq d^*(M_{p^k}),$$
but this contradicts (\ref{r2}) due to Propositions 2.2 and 3.5.

We conclude that our assumption is false, so $G$ is a Dedekind group.
\hfill\rule{1,5mm}{1,5mm}\\

We end this section with an extension of Theorem D to the class of nilpotent groups. This result holds due to the multiplicativity property of $d^*.$\\

\textbf{Corollary 3.7} \textit{Let $p_1<p_2<\ldots <p_k$ be prime numbers, $G$ be a nilpotent group and $P_1, P_2, \ldots, P_k$ be its Sylow subgroups such that $$|P_i|=p_i^{n_i}, \text{ where } n_i\geq 3, \ \forall \ i\in \{1, 2,\ldots, k\}.$$ 
\begin{itemize}
\item[i)] Assume that $p_1=2$ and $n_1=3$. If $d^*(P_1)>d^*(D_8)$ and $$d^*(P_i)>d^*(M_{p_i^{n_i}}), \ \forall \ i\in \{2,3,\ldots, k\},$$ then $G$ is a Dedekind group;
\item[ii)] Assume that $n_1\geq 4$ if $p_1=2$, and $n_1\geq 3$ if $p_1$ is odd. If 
 $$d^*(P_i)>d^*(M_{p_i^{n_i}}), \ \forall \ i \in\{ 1,2,\ldots, k\},$$ then $G$ is a Dedekind group.
\end{itemize}}

The lower bounds that appear in the statements of Theorem D and Corollary 3.7 are the best possible ones since $D_8$ and $M_{p^n}$ are non-Dedekind groups. Note that Corollary 3.7 does not hold for non-nilpotent supersolvable groups. More exactly, regarding item \textit{i)}, we can take $G\cong C_{27}\rtimes Q_8$ (SmallGroup(216, 4)). Its Sylow subgroups are Dedekind groups, but $G$ is non-Dedekind since $d'(G)=\frac{2}{11}$. In what concerns  item \textit{ii)}, we can choose $G_1\cong C_3^3\rtimes C_4^2$ (SmallGroup(432, 425)) and the Zassenhaus metacyclic group (see Theorem 11 of \cite{17}) $$G_2=\langle x, y \mid x^{6859}=y^{27}=1, yx=x^{956}y \rangle\cong C_{6859}\rtimes C_{27}.$$ The Sylow subgroups of both groups are Dedekind. However $d'(G_1)=\frac{89}{224}$ and $G_2$ has a non-normal subgroup $H\cong C_{27}$, so both groups are non-Dedekind.     

\section{Open problems}

We end our paper by enumerating some questions that can be further explored.\\

\textbf{Problem 4.1.} Let $a, b$ be integers such that $1\leq a<b$. In Section 2 we proved that any ratio of the form $\frac{a}{a+1}$ can be achieved by computing $d'$ for specific  non-nilpotent supersolvable groups having one conjugacy class of non-normal subgroups. We also showed that each ratio of the form $\frac{a}{b}$ can be achieved as the limit of $(d'(G_n))_{n\geq 1}$ where $(G_n)_{n\geq 1}$ is a sequence of nilpotent groups formed as direct products of some modular $p$-groups. Hence, we pose the following question:\\ 

\textit{Is it true that there is a group $G$ such that $d'(G)=\frac{a}{b}$?}\\

To get an affirmative answer, it suffices to search for $b-a$ groups $G_1, G_2, \ldots, G_{b-a}$ of coprime orders such that $$d'(G_i)=\frac{a+i-1}{a+i}, \ \forall \ i\in\{ 1, 2,\ldots, b-a\}.$$ By taking $G=\xmare{i=1}{b-a}G_i$ and making use of the multiplicativity of $d'$, we would get a positive answer.\\

\textbf{Problem 4.2.} Let $G$ be a group with $d^*(G)>\frac{2}{3}$. According to Theorem B, it follows that $G$ is a nilpotent group. Moreover, due to the same result, if $2\nmid |G|$, then $G$ is an Iwasawa group. Assume that $2\mid |G|$. Then $G$ has a Sylow 2-subgroup $P$ with $$d^*(P)\geq d^*(G)>\frac{2}{3}.$$

\textit{What can be said about the structure of $P$?}\\

Note that if $\frac{2}{3}$ is replaced with $\frac{4}{5}$, an answer is given by Iwasawa's result that is recalled after proving Theorem B. One can find a 2-group $P$ such that $d^*(P)\in (\frac{2}{3}, \frac{4}{5})$. For instance, by checking the groups of order 16, then $P\cong C_2^2\rtimes C_4$ (SmallGroup(16, 3)) and $d^*(P)=d'(P)=\frac{17}{23}$, or $P\cong C_2\times D_8$ (SmallGroup(16, 11)) and $d^*(P)=d'(P)=\frac{27}{35}$.\\

\textbf{Problem 4.3.} Let $n\geq 3$ be an integer. Theorem D indicates the second maximum value attained by $d^*$ on the class of $p$-groups of order $p^n$.\\ 

\textit{Which is the second maximum value achieved by $d'$ on the same class of groups?}\\

Assume that $p=2$ and $G$ is a non-Dedekind group of order $2^n$. If $n=3$, one can check that $d'(G)\leq \frac{4}{5}=d'(D_8)$, while if $n=4$, we have $d'(G)\leq\frac{10}{11}=d'(M_{16})$. Also, for $n\in\{5, 6, 7\}$, GAP computations show that $d'(G)\leq d'(C_2^{n-5}\times (D_8\circ Q_8))$. We mention that $C_2^{k}\times (D_8\circ Q_8)$ is a generalized extraspecial group of order $2^{n}$, where $n=5+k$ and $k\geq 0$.\\

\textit{Which are the minimum values achieved by  $d'$ and $d^*$ on the same class of groups?}\\

If $p=2$, it can be checked that $d^*(G)\geq d^*(D_{2^n})$ for any $n\in \{ 3, 4, 5, 6, 7\}$. Since any section of $D_{2^n}$ is abelian or isomorphic to $D_{2^k}$ where $k\leq n$, by using Proposition 2.2, we deduce that   $d'(D_{2^n})=d^*(D_{2^n})$. Hence, we also have that $d'(G)\geq d'(D_{2^n})$ for the same values of $n$.

\bigskip\noindent {\bf Acknowledgements.} The authors  are grateful to the reviewer for his/her remarks which improved the initial version of the paper. 

\bigskip\noindent {\bf Declarations}

\bigskip\noindent {\bf  Funding.} The authors did not receive support from any organization for the submitted work.

\bigskip\noindent {\bf Conflicts of  interests.} The authors declare that there is  no conflict of interest.

\bigskip\noindent {\bf  Data availability statement.} The manuscript has no associated data.

\vspace*{3ex}
\small
\begin{minipage}[t]{7cm}
Mihai-Silviu Lazorec \\
Faculty of  Mathematics \\
"Al. I. Cuza" University \\
Ia\c si, Romania \\
e-mail: {\tt silviu.lazorec@uaic.ro}
\end{minipage}
\hspace{2cm}
\begin{minipage}[t]{7cm}
Marius T\u arn\u auceanu \\
Faculty of  Mathematics \\
"Al. I. Cuza" University \\
Ia\c si, Romania \\
e-mail: {\tt tarnauc@uaic.ro}
\end{minipage}
\vspace{1cm}

\appendix
\section{GAP code}
\label{gap}

The following GAP code may be used to check most of our  results. 

\begin{lstlisting}[language=GAP, mathescape=true]
# A function that returns a list containing the sections of a group $G$

sections:=function(G)
    local u, H, N;
    u := [];
    for H in List(ConjugacyClassesSubgroups(G), Representative) do
        for N in NormalSubgroups(H) do
            Add(u, FactorGroup(H,N));
        od;
    od;
    return u;
end;

# A function that computes $d'(G)$

dG:=function(G)
    local c;
    c := ConjugacyClassesSubgroups(G);
    return Length(c) / Sum(c, Size);
end;

# A function that computes $d^*(G)$ 
# It should be applied to the output of sections(G), i.e. dStarG(sections(G))

dStarG:=function(l)
    return Minimum(List(l, dG));
end;


\end{lstlisting}
\end{document}